\documentclass[amstex,12pt,russian,amssymb]{article}

\usepackage{mathtext}
\usepackage[cp1251]{inputenc}
\usepackage[T2A]{fontenc}
\usepackage[russian]{babel}
\usepackage[dvips]{graphicx}
\usepackage{amsmath}
\usepackage{amssymb}
\usepackage{amsxtra}
\usepackage{latexsym}
\usepackage{ifthen}

\renewcommand{\baselinestretch}{1.3}

\begin{document}

\def\Xint#1{\mathchoice
   {\XXint\displaystyle\textstyle{#1}}%
   {\XXint\textstyle\scriptstyle{#1}}%
   {\XXint\scriptstyle\scriptscriptstyle{#1}}%
   {\XXint\scriptscriptstyle\scriptscriptstyle{#1}}%
   \!\int}
\def\XXint#1#2#3{{\setbox0=\hbox{$#1{#2#3}{\int}$}
     \vcenter{\hbox{$#2#3$}}\kern-.5\wd0}}
\def\dashint{\Xint-}

\renewcommand{\abstractname}{}

\noindent
\renewcommand{\baselinestretch}{1.1}
\normalsize \large \noindent

\newtheorem{theorem}{Теорема}[section]
\newtheorem{lemma}{Лемма}[section]
\newtheorem{proposition}{Предложение}[section]
\newtheorem{corollary}{Следствие}[section]
\newtheorem{definition}{Определение}[section]

\newtheorem{example}{Пример}[section]
\newtheorem{remark}{Замечание}[section]
\newcommand{\keywords}{\textbf{Ключевые слова и фразы. }\medskip}
\newcommand{\subjclass}{\textbf{MSC 2000. }\medskip}
\renewcommand{\abstract}{\textbf{Аннотация. }\medskip}
\numberwithin{equation}{section}

\title{Об одном свойстве кольцевых  $Q$-гомеоморфизмов относительно $p$-модуля}

\author{ Салимов Р.Р.}

\medskip

 {УДК 517.5}

{\bf  Салимов Р.Р.}  {(Ин-т прикладной математики и механики НАН
Украины, Донецк)}
\medskip

О конечной  липшицевости классов Орлича-Соболева

\medskip

On finite lipschitz Orlicz Sobolev classes

\begin{abstract}
Найдено   достаточное условие конечной
липшицевости гомеоморфизмов  класса Орлича-Соболева $W_{loc}^{1,\varphi}$ при условии типа Кальдерона на $\varphi$.

It is found a sufficient condition of finite Lipschitz of homeomorphisms of the Orlicz-Sobolev class
$W_{loc}^{1,\varphi}$ under a condition of the Calderon type.

\end{abstract}

\section{Введение}\label{int}
\vskip10pt

Напомним некоторые определения. Борелева функция
$\rho:\mathbb{R}^n\to[0,\infty]$ называется {\it допустимой} для
семейства кривых $D$ в $\mathbb{R}^n$, $n\geqslant 2$, пишут $\rho\in{\rm
adm}\,\Gamma$, если

\begin{equation}
\int\limits_{\gamma}\rho(x)\,ds\geqslant 1
\end{equation}
для всех $\gamma\in\Gamma$. Пусть $p\geqslant1$. Тогда
$p${\it-модулем} семейства кривых $\Gamma$ называется  величина

\begin{equation}
M_{p}(\Gamma)=\inf\limits_{\rho\in{\rm adm}\,\Gamma}
\int\limits_{\mathbb{R}^n}\rho^{p}(x)\,dm(x).
\end{equation}
Здесь
$m$ обозначает меру Лебега в $\mathbb{R}^n$.

Пусть $D$ -- область в  $\mathbb{R}^n$, $n\geqslant 2$. Предположим, что $n-1<p<n$ и \begin{equation}\label{eq1**}
M_p(f\Gamma)\leqslant K\,M_p(\Gamma)\end{equation} для произвольного
семейства $\Gamma$ кривых $\gamma$ в области $D$. При
предположении, что $f$ в (\ref{eq1**}) является гомеоморфизмом,
Герингом  было установлено, что отображение $f$ является {\it липшицевым}, другими словами, при некоторой постоянной $C>0$
и всех $x_0\in D$ справедлива оценка

\begin{equation}
\limsup\limits_{x\to x_0}\frac{|f(x)-f(x_0)|}{|x-x_0|}\leqslant
K^{\frac{1}{n-p}},
\end{equation}
см., напр., теорему 2 в \cite{Ge1}.

\medskip

Пусть  $D$ -- область в ${\Bbb R}^n$, $n\geqslant2$. Напомним, что гомеоморфизм $f:D\rightarrow {\Bbb R}^n$  называется {\it отображением с конечным
искажением,} если $f\in W^{1,1}_{\rm loc}$ и
\begin{equation}\label{eqOS1.3}\Vert
f^{\,\prime}(x)\Vert^n\le K(x)\cdot J_f(x)\end{equation} для
некоторой почти всюду конечной функции $K(x)\geqslant 1,$ где $f^{\,\prime}(x)$
якобиева матрица $f,$ $\Vert f^{\,\prime}(x)\Vert$ -- её операторная
норма: $\Vert f^{\,\prime}(x)\Vert=\sup\limits_{|h|=1} |f^{\,\prime}(x)\cdot
h|$ и  $J_f(x)=\det f^{\,\prime}(x)$ -- якобиан отображения $f$.

 Пусть $p\in (1, \infty)$. В дальнейшем, полагаем

\begin{equation}
\label{eq4.1.4} K_{p}(x,f)= \left \{\begin{array}{rr}
\frac{\Vert f'(x)\Vert^p}{J(x,f)}, & \text{\rm  если} \ J(x,f) \neq 0 \\
1, & \text{\rm  если} \ f^\prime(x)=0 \\ \infty, & {\rm } \text{в остальных точках}\,.
\end{array} \right. \end{equation}
Впервые
понятие отображения с конечным искажением введено в случае плоскости
для $f\in W^{1,2}_{\rm loc}$ в работе \cite{IS}, см. также
\cite{IM}.

Следуя Орличу, для заданной выпуклой
возрастающей функции $\varphi:[0,\infty)$ $\rightarrow[0,\infty)$,
$\varphi(0)=0$, обозначим символом $L^{\varphi}$ пространство всех
функций $f:D\rightarrow{\Bbb R},$ таких что
%
$$\int\limits_{D}\varphi\left(\frac{|f(x)|}
{\lambda}\right)\,dm(x)<\infty$$
при некотором
$\lambda>0,$ см., напр., \cite{KR}. Здесь  $m$ -- мера
Лебега в  ${\Bbb R}^n$. Пространство $L^{\varphi}$
называется {\it пространством Орлича}.

{\it Классом Орлича--Соболева} $W^{1,\varphi}_{\rm loc}(D)$
называется класс всех локально интегрируемых функций $f,$ заданных в
$D,$ с первыми обобщёнными производными по Соболеву, градиент
$\nabla f$ которых принадлежит классу Орлича локально в области $D.$
Если же, более того, $\nabla f$ принадлежит классу Орлича в области
$D,$ мы пишем $f\in W^{1,\varphi}(D).$ Заметим, что по определению
$W^{1,\varphi}_{\rm loc}\subset W^{1,1}_{\rm loc}.$ Как обычно, мы
пишем $f\in W^{1,p}_{\rm loc},$ если $\varphi(t)=t^p$, $p\geqslant 1.$
Известно, что непрерывная функция $f$ принадлежит классу
$W^{1,p}_{\rm loc}$ тогда и только тогда, когда $f\in ACL^{p}$,
т.е., если $f$ локально абсолютно непрерывна на почти всех прямых,
параллельных координатным осям, а первые частные производные $f$
локально интегрируемы в степени $p$ в области $D,$ см., напр.,
\cite[разд.~1.1.3]{Maz}.

Далее, если $f$ -- локально интегрируемая вектор--функция $n$
вещественных переменных $x_1,\ldots,x_n,$ $f=(f_1,\ldots,f_m),$
$f_i\in W_{\rm loc}^{1,1},$ $i=1,\ldots,m,$ и
%
$$\int\limits_{D}\varphi\left(|\nabla
f(x)|\right)\,dm(x)<\infty\,,$$
где $|\nabla
f(x)|=\sqrt{\sum\limits_{i=1}^m\sum\limits_{j=1}^n\left(\frac{\partial
f_i}{\partial x_j}\right)^2},$ то мы снова пишем $f\in
W^{1,\varphi}_{\rm loc}.$ Мы также используем обозначение
$W^{1,\varphi}_{\rm loc}$ в случае более общих функций $\varphi,$
чем в классах Орлича, всегда предполагающих выпуклость функции
$\varphi$ и ее нормировку $\varphi(0)=0.$

Отметим, что классы Орлича--Соболева сейчас, как и ранее, изучаются в самых различных аспектах
многими авторами, см., напр., \cite{ARS, AC,  Ca, Ci, Do, GM$_*$, Hs, IKO, KRSS1, Ko, KKM, KP,   LM, LL, On$_4$, Tu}
и \cite{Vui}.

\section{Свойства классов Орлича-Соболева }\label{2}

\begin{theorem}\label{thOS4.0}
Пусть $\Omega$ -- открытое множество в ${\Bbb
R}^n$, $n\geqslant3$, $f:\Omega\to{\Bbb R}^n$ -- непрерывное
открытое отображение класса $W^{1,\varphi}_{\rm loc}(\Omega)$, где
$\varphi:(0,\infty)\to(0,\infty)$ -- неубывающая функция,
удовлетворяющая условию
\begin{equation}\label{eqOS4.1}
\int\limits_{t_*}^{\infty}\left[\frac{t}{\varphi(t)}\right]^
{\frac{1}{n-2}}dt<\infty.\end{equation}
Тогда отображение $f$
имеет почти всюду полный дифференциал в $\Omega$.
\end{theorem}

\begin{remark}
В частности, заключение теоремы \ref{thOS4.0}
имеет место, если $f\in W^{1,p}_{\rm loc}$ при некотором $p>n-1$.
Последнее утверждение -- результат Вяйсяля, см. лемму 3 в
\cite{Va_65}. Теорема \ref{thOS4.0} является также распространением
в пространство хорошо известной теоремы Меньшова-Геринга-Лехто на
плоскости, см., напр., \cite{Me}, \cite{GL} и \cite{LV}.
\end{remark}

\begin{theorem}\label{thdiff}
Пусть $\Omega$ -- открытое множество в ${\Bbb
R}^n$, $n\geqslant3$, $f:\Omega\to{\Bbb R}^n$ -- непрерывное
открытое отображение класса $W^{1,\varphi}_{\rm loc}(\Omega)$, где
$\varphi:(0,\infty)\to(0,\infty)$ -- неубывающая функция,
удовлетворяющая условию
\begin{equation}\label{eqOS4.1}
\int\limits_{t_*}^{\infty}\left[\frac{t}{\varphi(t)}\right]^
{\frac{1}{n-2}}dt<\infty.\end{equation}
Тогда отображение $f$
имеет почти всюду полный дифференциал в $\Omega$.
\end{theorem}

\begin{theorem}\label{thN}
Пусть $U$ -- открытое множество в ${\Bbb R}^n$,
$n\geqslant3$, $\varphi:(0,\infty)\to(0,\infty)$ -- неубывающая
функция, такая что для некоторого $t_*\in (0, \infty)$
\begin{equation}\label{eqOS3.3} \int\limits_{t_*}^{\infty}\left[\frac{t}{\varphi(t)}\right]^
{\frac{1}{n-2}}dt<\infty.\end{equation} Тогда любое непрерывное
отображение $f:U\to{\Bbb R}^m$, $m\geqslant1$, класса
$W^{1,\varphi}_{\rm loc}$ обладает $(N)$-свойством, более того,
локально абсолютно непрерывно относительно $(n-1)$-мерной
хаусдорфовой меры на почти всех гиперплоскостях $\mathcal{P}$,
параллельных произвольной фиксированной гиперплоскости ${\mathcal
P}_0$. Кроме того, на почти всех таких $\mathcal{P}$,
$H^{n-1}(f(E))=0$, если $|\nabla f|=0$ на $E\subset\mathcal{P}$.
\end{theorem}

Заметим, что, если условие  вида (\ref{eqOS3.3}) имеет место для
некоторой неубывающей функции $\varphi$, то функция
$\varphi_c(t)=\varphi(c\,t)$ при $c>0$ также удовлетворяет
соотношению (\ref{eqOS3.3}). Кроме того, хаусдорфовы меры являются
квазиинвариантными при квазиизометриях.

\begin{corollary}\label{corOS3.2}{\it  При условии
(\ref{eqOS3.3}) любое непрерывное отображение $f\in
W^{1,\varphi}_{\rm loc}$ обладает $(N)$-свойством относительно
$(n-1)$-мерной меры Хаусдорфа, более того, локально абсолютно
непрерывно на почти всех сферах $S$ с центром в заданной
предписанной точке $x_0\in{\Bbb R}^n$. Кроме того, на почти всех
таких сферах $S$ выполнено условие $H^{n-1}(f(E))=0$ как только
$|\nabla f|=0$ на множестве $E\subset S$.}
\end{corollary}

\section{Модули семейств поверхностей }\label{Sal_5}
\vskip10pt

Следуя \cite[разд.~9.2, гл.~9]{MRSY}, далее {\it $k$-мерной
поверхностью} $S$ в ${\Bbb R}^n$ называется произвольное
неп\-ре\-рывное отображение $S:\omega\to{\Bbb R}^n$, где $\omega$ --
открытое множество в $\overline{{\Bbb R}^k}:={\Bbb
R}^k\cup\{\infty\}$ и $k=1,\ldots,n-1$. {\it Функцией кратности}
поверхности $S$ называется число прообразов
%
$$N(S,y)={\rm card}\,S^{\,-1}(y)={\rm card}\,\{x\in\omega:\
S(x)=y\},\quad y\in{\Bbb R}^n\,.$$
%
Другими словами, символ $N(S,y)$ обозначает кратность накрытия точки
$y$ поверхностью $S$. Известно, что функция кратности является
полунепрерывной снизу, и, значит, измерима относительно произвольной
хаусдорфовой меры $H^k,$ см., \cite[разд.~9.2]{MRSY}.

\medskip
Для борелевской функции $\rho:{\Bbb R}^n\to[0,\infty]$ ее {\it
интеграл над поверхностью} $S$ определяется равенством
\begin{equation}\label{admn1}
\int\limits_S\rho\,d{\mathcal A}:= \int\limits_{{\Bbb
R}^n}\rho(y)\:N(S,y)\,dH^ky\,.
\end{equation}

Пусть $\Gamma$ -- семейство $k$-мерных поверхностей $S$. Борелева
функция $\rho:{\Bbb R}^n\to[0,\infty]$ называется {\it допустимой}
для семейства $\Gamma$, пишут $\rho\in{\rm adm}\,\Gamma$, если
$$\int\limits_S\rho^k\,d{\mathcal A}\geqslant 1$$
для каждой поверхности $S\in\Gamma.$ Пусть $p\in(1,\infty)$ --
заданное фиксированное число. Тогда {\it $p$-модулем} семейства
$\Gamma$ называется величина
%
$$M_p(\Gamma)=\inf_{\rho\in{\rm adm}\,\Gamma}\int\limits_{{\Bbb
R}^n}\rho^p(x)\,dm(x)\,.$$

\medskip

\section{О емкости конденсатора }\label{Sal_5}
\vskip10pt

Следуя работе \cite{MRV}, пару $\mathcal{E}=(A,C)$, где $A\subset\mathbb{R}^n$
-- открытое множество и $C$ -- непустое компактное множество,
содержащееся в $A$, называем {\it конденсатором}. Конденсатор $\mathcal{E}$
называется  {\it кольцевым конденсатором}, если $G=A\setminus C$ --
кольцо, т.е., если $G$ -- область, дополнение которой
$\overline{\mathbb{R}^n}\setminus G$ состоит в точности из двух
компонент.  Говорят
также, что конденсатор $\mathcal{E}=(A,C)$ лежит в области $D$, если $A\subset
D$. Очевидно, что если $f:D\to\mathbb{R}^n$ -- непрерывное, открытое
отображение и $\mathcal{E}=(A,C)$ -- конденсатор в $D$, то $(fA,fC)$ также
конденсатор в $fD$. Далее $f\mathcal{E}=(fA,fC)$.

\medskip
Функция $u:A\to \mathbb{R}$ {\it абсолютно непрерывна на прямой}, имеющей непустое пересечение с $A$, если она абсолютно непрерывна на любом отрезке этой прямой, заключенном в $A$. Функция $u:A\to \mathbb{R}$ принадлежит классу ${\rm ACL}$ ({\it абсолютно непрерывна на почти всех прямых}), если она абсолютно непрерывна на почти всех прямых, параллельных любой координатной оси.

Обозначим через $C_0(A)$  множество
непрерывных функций $u:A\to\mathbb{R}^1$ с компактным носителем,
$W_0(\mathcal{E})=W_0(A,C)$ -- семейство неотрицательных функций
$u:A\to\mathbb{R}^1$ таких, что 1) $u\in C_0(A)$, 2)
$u(x)\geqslant1$ для $x\in C$ и 3) $u$ принадлежит классу ${\rm
ACL}$. Также обозначим

\begin{equation}
\vert\nabla
u\vert={\left(\sum\limits_{i=1}^n\,{\left(\frac{\partial u}{\partial x_i}\right)}^2
\right)}^{1/2}.
\end{equation}

При $p\geqslant1$ величину
\begin{equation}
{\rm cap_p}\,\mathcal{E}={\rm cap_p}\,(A,C)=\inf\limits_{u\in W_0(\mathcal{E})}\,
\int\limits_{A}\,\vert\nabla u\vert^p\,dm(x)
\end{equation}
называют
{\it $p$-ёмкостью} конденсатора $\mathcal{E}$. В дальнейшем при  $p>1$ мы
будем использовать равенство
\begin{equation}\label{EMC}
{\rm cap_p}\,\mathcal{E}=M_p(\Delta(\partial A,\partial C; A\setminus C)),\ \
 \end{equation} см. \cite{G}, \cite{H} и
\cite{Sh}.

Известно, что при $1\leqslant p<n$ \begin{equation}\label{maz} {\rm
cap_p}\,\mathcal{E}\geqslant n{\nu}^{\frac{p}{n}}_n
\left(\frac{n-p}{p-1}\right)^{p-1}\left[m(C)\right]^{\frac{n-p}{n}}\end{equation}
где ${\nu}_n$ -  объем  единичного шара  в ${\Bbb R}^n,\,\,$ см., напр., неравенство (8.9)
в  \cite{MazIS}.

При $n-1<p\leqslant n$ имеет место оценка
 \begin{equation}\label{krd}
\left({\rm cap_p}\,\,
\mathcal{E}\right)^{n-1}\,\geqslant\,\gamma\,\frac{d(C)^{p}}{m(A)^{1-n+p}}\,\,,
\end{equation}
где $d(C)$ - диаметр компакта $C$,  $\gamma$ - положительная константа, зависящая только от
размерности $n$ и $p\,,$ см. предложение 6  в \cite {Kru}.

 \section{Нижние $Q$-гомеоморфизмы относительно $p$-модуля}\label{8}
\vskip10pt

Говорят, см. \cite[разд.~9.2]{MRSY}, что измеримая по Лебегу функция
$\rho:{\Bbb R}^n \rightarrow [0,\infty]$ является {\it обобщённо
$p$-до\-пус\-ти\-мой} для семейства $\Gamma$, состоящего из
$(n-1)$ - мерных поверхностей $S$ в ${\Bbb R}^n$, пишут $\rho\in{\rm
ext}_p\,{\rm adm}\,\Gamma$, если
\begin{equation}
\int\limits_S\rho^{n-1}(x)\,d{\mathcal
A}\geqslant 1
\end{equation}
для $p$-поч\-ти всех $S\in\Gamma.$

В работе \cite{Ge}, разд. 13, Ф. Геринг определил $K$-ква\-зи\-кон\-фор\-м\-ное отображение как гомеоморфизм,
изменяющий модуль кольцевой области не более чем в $K$ раз. Следующее понятие мотивировано кольцевым определением Геринга.

Пусть $D$ и $D^{\prime}$ -- области в ${\Bbb R}^n$, $n\geqslant2$,
$x_0\in D$, $Q:D\to(0,\infty)$ измеримая по Лебегу функция.
Гомеоморфизм $f:D\to D^{\prime}$ будем называть {\it нижним
$Q$-го\-мео\-морфизмом относительно $p$-мо\-ду\-ля в точке}\index{нижний
$Q$-гомеоморфизмом в точке} $x_0,$ если
\begin{equation}\label{Sal_eqOS1.10} M_p\left(f\Sigma_{\varepsilon}\right)\geqslant\inf\limits_{\rho\in {\mathrm {ext_p\,adm}}\,
\Sigma_{\varepsilon}}\int\limits_{R}\frac{\rho^p(x)}{Q(x)}\ dm(x)\end{equation} для каждого кольца
$$R=R(x_0,\varepsilon, \varepsilon_0)=
\left\{x\in {\Bbb R}^n:\varepsilon<|x-x_0|<\varepsilon_0\right\},\varepsilon\in
(0,\varepsilon_0),\varepsilon_0\in(0,d_0),$$
где
$d_0=\mathrm {dist} (x_0, \partial D)\,,$ а $\Sigma_{\varepsilon}$ обозначает семейство всех  сфер
\begin{equation}
S(x_0,r)=\left\{x\in {\Bbb R}^n:|x-x_0|=r\right\}\,,\qquad r\in(\varepsilon,\varepsilon_0)\,.
\end{equation}

\medskip

\medskip

Следующий критерий нижних $Q$-гомеоморфизмов, см.  \cite[Теорема~6.1]{GoSa},  впервые был доказан при $p=n$ в работе \cite{KR$_1$},
теорема 2.1, см. также монографию \cite{MRSY}, теорема 9.2.

\bigskip

\begin{lemma}{}\label{SalLEMin1} {\it Пусть $D$ -- область в $\mathbb{R}^n$, $n\geqslant2$, $x_0\in\overline{D}$,
и пусть $Q:D\to(0,\infty)$ -- измеримая функция. Гомеоморфизм
$f:D\to\mathbb{R}^n$ является нижним $Q$-го\-мео\-морфизмом в точке
$x_0$ относительно $p$-мо\-ду\-ля при $p>n-1$ тогда и только тогда,
когда
\begin{equation}\label{Sal_eq2.1.13}M_p(f\Sigma_{\varepsilon})\geqslant\int\limits_{\varepsilon}^{\varepsilon_0}
\frac{dr}{||\,Q||\,_{\frac{n-1}{p-n+1}}(r)}\quad\forall\
\varepsilon\in(0,\varepsilon_0)\,,\quad\varepsilon_0\in(0,d_0)\,,\end{equation} где
$d_0=\sup\limits_{x\in
D}\,|x-x_0|\,,$ $\Sigma_{\varepsilon}$ -- семейство
всех пересечений сфер $S(x_0,r)=\{x\in\mathbb{R}^n:|\,x-x_0|=r\}$,
$r\in(\varepsilon,\varepsilon_0)$, с $D$, и
$$\|Q\|_{\frac{n-1}{p-n+1}}(r)=\left(\int\limits_{D(x_0,r)}
Q^{\frac{n-1}{p-n+1}}(x)\,d{\cal A}\right)^\frac{p-n+1}{n-1}\, ,$$ где $D(x_0,r)=\{x\in D:|\,x-x_0|=r\}=D\cap S(x_0,r)$.
Инфимум в (\ref{Sal_eqOS1.10}) достигается только для функции}
\begin{equation}\label{Sal_eq2.1.16}\rho_0(x)=\frac{Q(x)}{\|Q\|_{\frac{n-1}{p-n+1}}(|x-x_0|)}\,.\end{equation}
\end{lemma}

\medskip

\begin{lemma}\label{SalLEMin12} {\it Пусть $D$ -- область в $\mathbb{R}^n$, $n\geqslant2$, $x_0\in D$,
и пусть $Q:D\to(0,\infty)$ -- измеримая функция и
$f:D\to\mathbb{R}^n$ --  нижний  $Q$-го\-мео\-морфизмом в точке
$x_0$ относительно $p$-мо\-ду\-ля при $p>n-1$. Тогда  имеет место оценка
\begin{equation}\label{Sal_eq2.1.13}
M_{\frac{p}{p-n+1}}\left(\Delta(fS_{1}, fS_{2},
fD)\right)\ \leqslant \left(\int\limits_{r_1}^{r_2}
\frac{dr}{\Vert\,Q\Vert\,_{\frac{n-1}{p-n+1}}(r)}\right)^{-\frac{n-1}{p-n+1}}\,,\end{equation}
где  $S_j=S(x_0, r_j), j=1,2$.}
\end{lemma}

{\bf Доказательство.} Действительно, пусть $0<r_1<r_2<d(x_0,\partial D)$ и $S_i=S(x_0,r_i),$ $i=1,2.$ Согласно неравенствам Хессе и Цимера
(см., напр., \cite{Hes} и \cite{Zi}.),
\begin{equation}\label{eqOS6.1aa} M_{\frac{p}{p-n+1}}\left(f\left(\Delta(S_{1}, S_{2}, D)\right)\right)\leqslant
\frac{1}{M_p^{\frac{n-1}{p-n+1}}(f\left(\Sigma\right))}\,,\end{equation}
поскольку
$f\left(\Sigma\right)\subset\Sigma\left(f(S_1),f(S_2),f(D)\right),$
где $\Sigma$ обозначает совокупность всех сфер с центром в точке
$x_0,$ расположенных между сферами $S_1$ и $S_2,$ а
$\Sigma\left(f(S_1),f(S_2),f(D)\right)$ состоит из всех
$(n-1)$-мерных поверхностей в $f(D),$ отделяющих $f(S_1)$ и
$f(S_2).$ Из соотношения (\ref{eqOS6.1aa}) по лемме
\ref{SalLEMin1} вытекает  заключение леммы  \ref{SalLEMin12}.

\bigskip

\bigskip

 \section{Взаимосвязь нижних $Q$-гомеоморфизмов  с  классами  Орлича-Соболева}\label{8}
\vskip10pt

Напомним, что отображение $g:X\to Y$ между метрическими
пространствами  $X$ и $Y$ называется {\it липшицевым,} если ${\rm
dist}\,\left(g(x_1),g(x_2)\right)\leqslant M\cdot{\rm
dist}\,\left(x_1,x_2\right)$ для некоторой постоянной $M<\infty$ и
всех $x_1$, $x_2\in X$. Говорят, что отображение $g:X\to Y$ {\it
билипшицево}, если, оно, во-первых, липшицево, во-вторых,
$M^*\cdot{\rm dist}\,\left(x_1,x_2\right)\leqslant{\rm
dist}\,\left(g(x_1),g(x_2)\right)$ для некоторой постоянной
$M^*>0$ и всех $x_1$, $x_2\in X$.\medskip

Следующее утверждение является ключевым для дальнейшего
исследования.

\begin{theorem}\label{thOS4.1} {\it Пусть $D$ и $D'$ -- области в ${\Bbb
R}^n$, $n\geqslant3$, $\varphi:(0,\infty)\to(0,\infty)$ --
неубывающая функция, такая что при $t_*\in (0, \infty)$
\begin{equation}\label{eqOS4.1}
\int\limits_{t_*}^{\infty}\left[\frac{t}{\varphi(t)}\right]^
{\frac{1}{n-2}}dt<\infty.\end{equation} Тогда любой гомеоморфизм
$f:D\to D'$ конечного искажения класса $W^{1,\varphi}_{\rm loc}$
является нижним $K_{p}(x,f)$-гомеоморфизмом относительно $p$-модуля с $p>n-1$.} \end{theorem}

{\bf Доказательство.} Обозначим через $B$ (борелево) множество всех точек $x\in D,$
где отображение $f$ имеет полный дифференциал и $J_f(x)=\mbox{det}\,
f^{\,\prime}(x)\ne 0.$ Заметим, что множество $B$ представляет собой
не более чем счётное объединение борелевских множеств $B_l,$
$l=1,2,\ldots\,$, таких что отображения $f_l=f|_{B_l}$ являются
билипшицевыми гомеоморфизмами, см., напр., \cite[лемма~3.2.2]{Fe}.
Без ограничения общности, можно считать, что множества $B_l$ попарно
не пересекаются. Обозначим также через $B_*$ оставшееся множество
всех точек $x\in D,$ где $f$ имеет полный дифференциал, однако,
$f^{\,\prime}(x)=0.$

По теореме \ref{thOS4.0} множество $B_0:=D\setminus \left(B\bigcup
B_*\right)$ имеет меру Лебега нуль. Следовательно, по
\cite[теорема~9.1]{MRSY} имеем, что
 $H^{n-1}(B_0\cap S_r)=0$ для $p$-почти всех сфер $S_r:=S(x_0,r)$ с
центром в произвольной точке $x_0\in\overline{D},$ где "$p$-почти
всех" \ определяется в смысле $p$-модуля семейства поверхностей. Тогда,
в силу \cite[лемма~9.1]{MRSY}, $H^{n-1}(B_0\cap S_r)=0$ для почти
всех $r\in {\Bbb R}$ и по следствию \ref{corOS3.2} получаем, что
$H^{n-1}(f(B_0)\cap S^*_r)=0$ и $H^{n-1}(f(B_*)\cap S^*_r)=0$ для
почти всех $r\in {\Bbb R}$, где $S^*_r=f(S_r).$

\medskip
Заметим, что также $H^{n-1}(f(B_0)\cap S^*_r)=0$ и
$H^{n-1}(f(B_*)\cap S^*_r)=0$ для почти всех сфер $S_r:=S(x_0,r)$ в
смысле $p$-модуля семейства поверхностей. Действительно, пусть
$\Gamma_0$ -- подсемейство всех сфер $S_r:=S(x_0,r),$ для которых
либо $H^{n-1}(f(B_0)\cap S^*_r)>0,$ либо $H^{n-1}(f(B_*)\cap
S^*_r)>0.$ Обозначим через $R$ множество всех $r\in {\Bbb R},$ для
которых либо $H^{n-1}(f(B_0)\cap S^*_r)>0,$ либо $H^{n-1}(f(B_*)\cap
S^*_r)>0.$ В силу сказанного выше, $m_1(R)=0.$ Тогда по теореме
Фубини $m(E)=0,$ где $E=\{x\in D: |x-x_0|=r\in R\}.$ Функция
$\rho_1:{\Bbb R}^n\rightarrow [0, \infty],$ определённая символом
$\infty$ при $x\in E$ и равная нулю на оставшемся множестве
обобщенно $p$-допустима для семейства $\Gamma_0$. Таким образом, по
(9.18) в \cite{MRSY} $M_p(\Gamma_0)\leqslant \int\limits_E \rho_1^p
dm(x)=0,$ т.е., действительно, $M_p(\Gamma_0)=0$.

\medskip
По теореме Кирсбрауна, см. \cite[теорема~2.10.43]{Fe}, каждое
отображение $f_l$ может быть продолжено до липшицевского отображения
$\widetilde{f_l}:{\Bbb R}^n\rightarrow {\Bbb R}^n,$ которое по
теореме Радемахера--Степанова $\widetilde{f_l}$ дифференцируемо
почти всюду в ${\Bbb R}^n,$ см. \cite[теорема~3.1.6]{Fe}. В силу
единственности аппроксимативного дифференциала (см.
\cite[пункт~3.1.2]{Fe}), можно считать, что при всех $x\in B_l$
выполнено равенство $\widetilde{f_l}^{\,\prime}(x)=f^{\,\prime}(x).$

Пусть $\Gamma$ обозначает семейство всех пересечений сфер  $S_r$,
$r\in(\varepsilon,\varepsilon_0)$,
$\varepsilon_0<d_0=\sup\limits_{x\in D}\,|x-x_0|,$ с областью $D.$
Для произвольной функции $\rho_*\in{\rm adm}\,f(\Gamma),$ такой что
$\rho_*\equiv0$ вне $f(D)$, полагаем $\rho\equiv 0$ вне $D$ и на
$B_0,$ и
$$\rho(x)\ \colon=\ \rho_*(f(x))\Vert f^{\,\prime}(x)\Vert \qquad\text{при}\ x\in D\setminus
B_0=B\cup B_*$$  Рассуждая покусочно на каждом $B_l$,
$l=1,2,\ldots$, согласно  \cite[разд.~1.7.6]{Fe},  а также используя
геометрический смысл величины $\Vert
f^{\,\prime}(x)\Vert$ и её связь с 
якобианом отображения, см., напр., соотношения (2.5) и (2.6) гл. I
$\S\, 2$ в \cite{Re}, имеем, что
$$\int\limits_{S_r}\rho^{n-1}\,d{\cal A}=
\int\limits_{S_r}\rho_*^{n-1}(f(x))\Vert
f^{\,\prime}(x)\Vert^{n-1}\,d{\cal A}=$$
$$=\int\limits_{S_r}\rho_*^{n-1}(f(x))\cdot\frac{\Vert
f^{\,\prime}(x)\Vert^{n-1}}{\frac{d{\cal A_*}}{d{\cal A}}}\cdot
\frac{d{\cal A_*}}{d{\cal A}}\,d{\cal A}\geqslant
\int\limits_{S_r}\rho_*^{n-1}(f(x))\cdot \frac{d{\cal A_*}}{d{\cal
A}}\,d{\cal A}=$$
$$=\int\limits_{S^*_r}\rho^{n-1}_*(y)d{\cal A_*}\geqslant 1$$ для почти всех $S_r,$ и, следовательно,
$\rho\in{\mathrm{ext_p\,adm}}\,\Gamma.$ Используя замену переменных на
каждом $B_l$, $l=1,2,\ldots,$ см., напр., \cite[теорема~3.2.5]{Fe},
ввиду счётной аддитивности интеграла, получаем также оценку
$$\int\limits_{D}\frac{\rho^p(x)}{K_{p}(x,f)}\,dm(x)\ \leqslant\
\int\limits_{f(D)}\rho^p_*(x)\, dm(x)\,,$$ что и завершает
доказательство.

\begin{corollary}\label{corOS4.1} {\it Любой гомеоморфизм с конечным искажением в ${\Bbb R}^n$,
$n\geqslant3$, класса $W^{1,\alpha}_{\rm loc}$ при $\alpha>n-1$ является
нижним $K_{p}(x,f)$-гомеоморфизмом с $p>n-1$.} \end{corollary}

\subsection{ Конечная липшицевость классов Орлича-Соболева}\label{Sal_6}

Для непрерывного  отображения $f:D\to\mathbb{R}^n$ и
$x\in D\subseteq\mathbb{R}^n$, положим

\begin{equation}
L(x,f)=\limsup_{y\to
x}\frac{|f(y)-f(x)|}{|y-x|}\,.
\end{equation}
Говорят, что отображение
$f$ является {\it конечно липшицевым}, если
$$L(x,f)<\infty$$ для всех $x\in D$.

\begin{theorem}\label{thOS4.1} {\it Пусть $D$ и $D'$ -- области в ${\Bbb
R}^n$, $n\geqslant3$. Предположим, что $f: D\to D'$ -- гомеоморфизм
с конечным  искажением  класса $W^{1,\varphi}_{\rm loc}$  с условием (\ref{eqOS4.1}) и, кроме того,
при  $p\in \left(n, n+\frac{1}{n-2}\right)$

\begin{equation}
 k_p(x_0)=\limsup\limits_{\varepsilon\to 0}
\left(\dashint_{B(x_0,\varepsilon) }  [K_{p}(x,f)]^{\frac{n-1}{p-n+1}}\,  dm(x)\right)^{\frac{p-n+1}{n-1}}<\infty.
\end{equation}
Тогда

\begin{equation}
L(x_0,f)=\limsup_{x\to
x_0}\frac{|f(x)-f(x_0)|}{|x-x_0|}\leqslant c_{n,p} \cdot  k^{\gamma}_p(x_0)<\infty\,,\end{equation}
где $\gamma=\frac{n-1}{n(p-n+1)-p}$  и   $c_{n,p}$  -- положительная константа, зависящая только от размерности
пространства   $n$ и $p$.
} \end{theorem}

{\it Доказательство.} Рассмотрим сферическое кольцо
$R=R(x_0,\varepsilon_1, \varepsilon_2)$ с $0<\varepsilon_1<\varepsilon_2$ такое,
что $R(x_0,\varepsilon_1, \varepsilon_2)\subset D$. Тогда $\mathcal{E}=\left(B\left(x_0,\varepsilon_2\right),\overline{B\left(x_0,\varepsilon_1\right)}\right)$
-- кольцевой конденсатор в   $D$ и $f\mathcal{E}=\left(fB\left(x_0,\varepsilon_2\right),\overline{fB\left(x_0,\varepsilon_1\right)}\right)$
-- кольцевой конденсатор в   $D'$.

Пусть $\Gamma^{*}=\Delta(fS_{1}, fS_{2},
fR)$, где  $S_j=S(x_0, r_j), j=1,2$. Тогда согласно (\ref{EMC}), имеем
равенство

\begin{equation}\label{eqCM}
{\rm cap_{\frac{p}{p-n+1}}}\ f\mathcal{E}=M_{\frac{p}{p-n+1}}\left( \Gamma^{*}\right)\,.
\end{equation}

По лемме \ref{SalLEMin12} получаем, что
\begin{equation}\label{eqCM1}
{\rm cap_{\frac{p}{p-n+1}}}\ f\mathcal{E} \leqslant \left(\int\limits_{\varepsilon_1}^{\varepsilon_2}
\frac{dr}{\Vert\,K_{p}(x,f)\Vert\,_{\frac{n-1}{p-n+1}}(r)}\right)^{-\frac{n-1}{p-n+1}}\,,
\end{equation}
где $\Vert\,K_{p}(x,f)\Vert\,_{\frac{n-1}{p-n+1}}(r)=
\left(\int\limits_{S(x_0,r)} [K_{p}(x,f)]^{\frac{n-1}{p-n+1}}\,d{\cal A}\right)^\frac{p-n+1}{n-1}\,.$

Заметим, что

\begin{equation}
\varepsilon_2-\varepsilon_1=\int\limits_{\varepsilon_1}^{\varepsilon_2} \Vert\,K_{p}(x,f)\Vert\,^{\frac {n-1}{p}}_{\frac{n-1}{p-n+1}}(r)\cdot \frac{dr}{ \Vert\,K_{p}(x,f)\Vert\,^{\frac {n-1}{p}}_{\frac{n-1}{p-n+1}}(r)}\,.
\end{equation}

И применяя неравенство Гельдера  с $q=\frac{p}{p-n+1}$  и $q'=\frac{p}{n-1}$ имеем

\begin{equation}\label{eq100fh1}
 \left(\int\limits_{\varepsilon_1}^{\varepsilon_2}
\frac{dr}{\Vert\,K_p(x,f)\Vert\,_{\frac{n-1}{p-n+1}}(r)}\right)^{-\frac{n-1}{p-n+1}}
\leqslant\frac{1}{(\varepsilon_2-\varepsilon_1)^{\frac{p}{p-n+1}}} \int\limits_{R}[K_p(x,f)]^{\frac{n-1}{p-n+1}} \,dm(x)\,.
\end{equation}

Комбинируя  неравенства (\ref{eq100fh1}) и (\ref{eqCM1}), получим

\begin{equation}\label{eq100fh}
 {\rm cap_{\frac{p}{p-n+1}}}\ f\mathcal{E} \leqslant
\frac{1}{(\varepsilon_2-\varepsilon_1)^{\frac{p}{p-n+1}}} \int\limits_{R}[K_p(x,f)]^{\frac{n-1}{p-n+1}} \,dm(x)\,.
\end{equation}

Далее, выбирая
$\varepsilon_1=2\varepsilon$ и $\varepsilon_2=4\varepsilon$, получим

\begin{equation}\label{eq101}{\rm cap_{\frac{p}{p-n+1}}}\ (fB(x_0,4\varepsilon),f\overline{B(x_0,2\varepsilon)})\leqslant\,
\frac{1}{(2\varepsilon)^{\frac{p}{p-n+1}}}\int\limits_{B(x_0,4\varepsilon)}[K_p(x,f)]^{\frac{n-1}{p-n+1}}\,dm(x)\,.
\end{equation}
С другой стороны,  в силу  неравенства (\ref{maz}) вытекает оценка

\begin{equation}\label{eq102} {\rm cap_{\frac{p}{p-n+1}}}\ (fB(x_0,4\varepsilon),f\overline{B(x_0,2\varepsilon)})
\geqslant c_{1}\left[m(fB(x_0,2\varepsilon))\right]^{\frac{n(p-n+1)-p}{n(p-n+1)}}\,,
\end{equation}
где    $c_{1}$  -- положительная константа, зависящая только от размерности
пространства   $n$ и $p.$

Комбинируя   (\ref{eq101}) и (\ref{eq102}), получаем, что
\begin{equation}\label{eq4.2} \frac{m(fB(x_0,2\varepsilon))}{m(B(x_0,2\varepsilon))}\leqslant c_{2}\,\left[
\dashint_{B(x_0,4\varepsilon)}\, [K_p(x,f)]^{\frac{n-1}{p-n+1}}\, dm(x)
\right]^{\frac{n(p-n+1)}{n(p-n+1)-p}}\,,\end{equation} где $c_{2}$ -
положительная постоянная зависящая только от $n$ и $p$.

Далее, выбирая в (\ref{eq100fh}) $\varepsilon_1=\varepsilon$ и
$\varepsilon_2=2\varepsilon$, получим
\begin{equation}\label{eq91}{\rm cap_{\frac{p}{p-n+1}}}\ (fB(x_0,2\varepsilon),f\overline{B(x_0,\varepsilon)})\leqslant\,
\frac{1}{\varepsilon^{\frac{p}{p-n+1}}}\int\limits_{B(x_0,2\varepsilon)}[K_p(x,f)]^{\frac{n-1}{p-n+1}}\,dm(x)\,.
\end{equation}
С другой стороны, в силу неравенства (\ref{krd}), получаем
\begin{equation}\label{eq10*!} {\rm cap_{\frac{p}{p-n+1}}}\ (fB(x_0,2\varepsilon),f\overline{B(x_0,\varepsilon)})
\geqslant
\left(c_3\frac{d^{\frac{p}{p-n+1}}(fB(x_0,\varepsilon))}{m^{1-n+\frac{p}{p-n+1}}(fB(x_0,2\varepsilon))}\right)^{\frac{1}{n-1}}\,,
\end{equation}
где    $c_3$  --   положительная константа, зависящая только от  $n$ и $p.$

Комбинируя   (\ref{eq91}) и (\ref{eq10*!}), получаем, что
$$
\frac{d(fB(x_0,\varepsilon))}{\varepsilon}\leqslant c_4
\left(\frac{m(fB(x_0,2\varepsilon))}{m(B(x_0,2\varepsilon))}
\right)^{i_1}\times
$$
\begin{equation}
\times \left( \dashint_{B(x_0,2\varepsilon)}
[K_p(x,f)]^{\frac{n-1}{p-n+1}}\, dm(x)\right)^{i_2}\,,
\end{equation}
где

$$
i_1=\frac{(1-n)(p-n+1)+p}{p}, \ \ \ \ i_2=\frac{(n-1)(p-n+1)}{p}
$$
и $c_4$  --   положительная константа, зависящая только от  $n$ и $p.$

Эта оценка вместе с  (\ref{eq4.2}) дает неравенство
$$
\frac{d(fB(x_0,\varepsilon))}{\varepsilon}\leqslant c_5\left(\dashint_{B(x_0,4\varepsilon)} [K_{p}(x,f)]^{\frac{n-1}{p-n+1}}
dm(x)\right)^{j_1}\times
$$

\begin{equation}\times \left(\dashint_{B(x_0,2\varepsilon)}\, [K_{p}(x,f)]^{\frac{n-1}{p-n+1}} dm(x)
\right)^{j_2}\,,
\end{equation}
где
$$
j_1=\frac{n\left((1-n)(p-n+1)+p\right)\left(p-n+1\right)}{p\left(n(p-n+1)-p\right)}, \ \ \  j_2=\frac{(n-1)(p-n+1)}{p}
$$
и   $c_5$  --   положительная константа, зависящая только от  $n$ и $p.$

Переходя к верхнему пределу при
$\varepsilon\to 0$, получаем
$$L(x_0,f)=\limsup\limits_{x\to
x_0}\frac{|f(x)-f(x_0)|}{|x-x_0|}\leqslant\limsup\limits_{\varepsilon\to
0}\frac{d(fB(x_0,\varepsilon))}{\varepsilon}\leqslant c\cdot
[k_p(x_0)]^{\frac{n-1}{n(p-n+1)-p}},
$$
где  $c$ - положительная  постоянная, зависящая только
от $n$ и $p$.

\bigskip

\begin{corollary}\label{thOS4.1} {\it Пусть $D$ и $D'$ -- области в ${\Bbb
R}^n$, $n\geqslant3$. Предположим, что $f: D\to D'$ -- гомеоморфизм
с конечным  искажением  класса $W^{1,\varphi}_{\rm loc}$  с условием (\ref{eqOS4.1}) и, кроме того,
при  $p\in \left(n, n+\frac{1}{n-2}\right)$
\begin{equation}
 \limsup\limits_{\varepsilon\to 0}
\dashint_{B(x_0,\varepsilon) }  [K_{p}(x,f)]^{\frac{n-1}{p-n+1}}\,  dm(x)<\infty\,  \ \ \ \ \ \ \ \forall x_0\in D.
\end{equation}
Тогда
гомеоморфизм  $f$ является конечно липшицевым.}
\end{corollary}

{\bf Замечание.} В соответствии с леммой 10.6 в \cite{MRSY} конечно
липшицевые отображения обладают $N$-свойством относительно
хаусдорфовых мер и, таким образом, являются абсолютно непрерывными
на кривых и поверхностях.

Построим пример  гомеоморфизма  с конечным искажением, не являющегося конечно липшицевым.

{\sl Пример.} Предположим, что $p\in \left(n, n+\frac{1}{n-2}\right)$. Пусть  $f:\mathbb{B}^n \to
\mathbb{B}^n$, где
$$f(x)=\frac{x}{|x|}\left(1+(p-n)\int\limits_{|x|}^{1}\frac{dt}{t^{p-n+1}\ln^{\frac{p-n+1}{n-1}}(\frac{e}{t})}\right)^{-\frac{1}{p-n}}$$
при  $x\neq 0$ и $f(0)=0$.

Касательная и радиальная дилатации $f$ на сфере $|x|=r$, $r\in (0,1)$, легко вычисляются:

$$ \delta_T=\frac{|f(x)|}{|x|}=\frac{\left(1+(p-n)\int\limits_{r}^{1}\frac{dt}{t^{p-n+1}\ln^{\frac{p-n+1}{n-1}}(\frac{e}{t})}\right)^{-\frac{1}{p-n}}}{r}$$
и

$$ \delta_r=\frac{\left(1+(p-n)\int\limits_{r}^{1}\frac{dt}{t^{p-n+1}\ln^{\frac{p-n+1}{n-1}}(\frac{e}{t})}\right)^{-\frac{p-n+1}{p-n}}}{r^{p-n+1}\ln^{\frac{p-n+1}{n-1}}(\frac{e}{t})}\,.$$

Заметим, что $\delta_T \geqslant \delta_r $ и
$$
\delta_T^{p-n+1}=\delta_r \ln^{\frac{p-n+1}{n-1}}(\frac{e}{r})\,.
$$
Следовательно, ввиду сферической симметрии мы видим, что

$$
K_{p}(x,f)=\frac{\delta_T^{p}}{\delta_T^{n-1}\delta_r}=\frac{\delta_T^{p-n+1}}{\delta_r}=\ln^{\frac{p-n+1}{n-1}}(\frac{e}{|x|})\,.
$$

Заметим, что

\begin{equation}
 \limsup\limits_{\varepsilon\to 0}
\dashint_{B(x_0,\varepsilon) }  [K_{p}(x,f)]^{\frac{n-1}{p-n+1}}\,  dm(x)=\infty\,.
\end{equation}

Тем не менее, как легко проверить по
правилу  Лопиталя, $\frac{|f(x)|}{|x|}\to \infty$ при  $x\to 0$,
т.е. гомеоморфизм $f$ не является липшицевым в нуле.

\medskip

Салимов Руслан Радикович

Институт прикладной математики и механики НАН Украины

ул. Розы  Люксембург 74, Донецк, 83114.

Рабочий телефон: 311-01-45

Email: salimov07@rambler.ru, ruslan623@yandex.ru,


\bigskip

\begin{thebibliography}{30}

\bibitem{Ge1} {\it Gehring F.W.} Lipschitz mappings and the $p$-capacity of ring in
$n$-space // Advances in the theory of Riemann surfaces (Proc. Conf.
Stonybrook, N.Y., 1969), Ann. of Math. Studies. -- 1971. --
\textbf{66}. -- P. 175–193.




\bibitem{IS} {\it Iwaniec T., Sver\'ak V.} On mappings with
integrable dilatation // Proc. Amer. Math. Soc. -- 1993. --
\textbf{118}. -- P. 181--188.


\bibitem{IM} {\it Iwaniec T., Martin G.} Geometrical Function Theory and
Non-Linear Analysis. -- Clarendon Press, Oxford,  2001.




\bibitem{KR} {\it Красносельский М.А., Рутицкий Я.Б.} Выпуклые функции и пространства
Орлича. -- Гос. издат. физ.-мат. лит., Москва, 1958.



\bibitem{Maz} {\it Мазья В.Г.} Пространства С.Л. Соболева. -- ЛГУ,
Ленинград, 1985. -- 416 с.














\bibitem{ARS}{\it Афанасьева Е.С.,  Рязанов В.И.,  Салимов Р.Р.} Об отображениях
в классах Орлича-–Соболева на римановых многообразиях // Укр. матем.
вісник. -- 2011. -- \textbf{8}, № 3. -- С. 319–-342.



\bibitem{AC}
{\it Alberico A., Cianchi A.} Differentiability properties of
Orlicz-Sobolev functions // Ark. Mat. -- 2005. -- \textbf{43}. -- P. 1--28.


\bibitem{Ca} {\it Calderon A.P.} On the differentiability of absolutely
continuous functions // Riv. Math. Univ. Parma. -- 1951. --
\textbf{2}. -- C. 203--213.



\bibitem{Ci}{\it Cianchi A.} A sharp embedding theorem for Orlicz-Sobolev
spaces // Indiana Univ. Math. J. -- 1996. -- \textbf{45}, no. 1. -- P. 39--65.

\bibitem{Do} {\it Donaldson T.} Nonlinear elliptic boundary-value problems in Orlicz-Sobolev
spaces // J. Diff. Eq. -- 1971. -- \textbf{10}. -- P. 507--528.








\bibitem{GM$_*$}
{\it Gossez J.-P., Mustonen V.} Variational inequalities in Orlicz-Sobolev spaces // Nonlinear Anal. Theory Meth.
Appl. -- 1987. -- \textbf{11}. --  \linebreak P. 379--392.

\bibitem{Hs}
{\it Hsini M.} Existence of solutions to a semilinear elliptic system through generalized
Orlicz-Sobolev spaces // J. Partial Differ. Equ. -- 2010. -- \textbf{23}, no. 2. -- P. 168--193.


\bibitem{IKO}
{\it Iwaniec T., Koskela P., Onninen J.} Mappings of
finite distortion: Compactness // Ann. Acad. Sci. Fenn. Ser. A1. Math. -- 2002. -- \textbf{27}, no. 2. -- P.
391--417.



\bibitem{KRSS1}{\it Kovtonyuk D.A., Ryazanov V.I., Salimov R.R., Sevost'yanov E.A.
} On mappings in the Orlicz-Sobolev classes // Ann. Univ. Bucharest,
Ser. Math. - 2012. - \textbf{3(LXI)}, no. 1. - P. 67-78.


\bibitem{Ko}
{\it Koronel J.D.} Continuity and $k$-th order differentiability in
Orlicz-Sobolev spaces: $W^kL_A$$"$ // Israel J. Math. -- 1976. -- \textbf{24}, no. 2. -- P. 119--138.

\bibitem{KKM}
{\it Kauhanen J., Koskela P., Maly J.} On functions with derivatives in a Lorentz
space // Manuscripta Math. -- 1999. -- \textbf{10}. -- P. 87--101.




\bibitem{KP}
{\it Khruslov E.Ya., Pankratov L.S.} Homogenization of the Dirichlet variational
problems in Sobolev-Orlicz spaces. -- Operator theory and its applications (Winuipeg, MB, 1998), 345-366, Fields
Inst. Commun., 25, Amer. Math. Soc,. Providence, RI, 2000.




\bibitem{LM} {\it Landes R., Mustonen V.} Pseudo-monotone mappings in Sobolev-Orlicz spaces
and nonlinear boundary value problems on unbounded domains // J.
Math. Anal. Appl. -- 1982. -- \textbf{88}. -- P. 25--36.


\bibitem{LL} {\it Lappalainen V. and Lehtonen A.} Embedding of Orlicz-Sobolev spaces in H\"older
spaces // Ann. Acad. Sci. Fenn. Ser. A1. Math. -- 1989. --
\textbf{14}, no. 1. -- P. 41--46.




\bibitem{On$_4$}
{\it Onninen J.} Differentiability of monotone Sobolev functions // Real. Anal. Exchange. -- 2000/2001. --
\textbf{26}, no. 2. -- P. 761--772.



\bibitem{Tu} {\it Tuominen H.} Characterization of Orlicz-Sobolev
space // Ark. Mat. -- 2007. -- \textbf{45}, no.~1. -- P. 123--139.



\bibitem{Vui}
{\it Vuillermot P.A.}
H\"older-regularity for the solutions of strongly nonlinear eigenvalue problems on
Orlicz-Sobolev space // Houston J. Math. -- 1987. -- \textbf{13}. -- P. 281--287.











\bibitem{Va_65}
{\it V\"{a}is\"{a}l\"{a} J.} Two new characterizations for quasiconformality // Ann. Acad. Sci. Fenn. Ser. A1
Math. -- 1965. -- \textbf{362.} -- P. 1--12.


\bibitem{Me} {\it Menchoff D.} Sur les differencelles totales des
fonctions univalentes // Math. Ann. -- 1931. -- \textbf{105}. -- P.
75--85.




\bibitem{GL} {\it Gehring F.W., Lehto O.} On the total
differentiability of functions of a complex variable // Ann. Acad.
Sci. Fenn. Ser. A1. Math. -- 1959. -- \textbf{272}. -- P. 3--8.


\bibitem{LV} {\it Lehto O., Virtanen K.} Quasiconformal Mappings in the
Plane. -- Springer--Verlag, New York, 1973.



\bibitem{MRSY} {\it Martio O., Ryazanov V., Srebro U., Yakubov E.} Moduli in Modern
Mapping Theory, Springer Monographs in Mathematics. -- Springer, New
York etc., 2009. -- 367 p.



\bibitem{MRV} {\it Martio O., Rickman S., and V\"{a}is\"{a}l\"{a} J.} Definitions for
quasiregular mappings // Ann. Acad. Sci. Fenn. Ser. A1. Math. --
1969. -- \textbf{448}. -- P. 1--40.





\bibitem{G} {\it Gehring F.W.} Quasiconformal mappings in Complex Analysis
and its Applications, V. 2, International Atomic Energy Agency,
Vienna, 1976.





\bibitem{H} {\it Hesse J.}  A $p$-extremal length and $p$-capacity equality //
Arc. Mat. -- 1975. -- \textbf{13}. -- P. 131-144.


\bibitem{Sh} {\it Shlyk V.A.} О равенстве $p$-емкости и $p$-модуля
// Сиб. матем. ж. -- 1993. -- \textbf{34}, №~6. -- С. 216--221.



\bibitem{MazIS}
V.~Maz'ya, {\em Lectures on isoperimetric and isocapacitary inequalities in the theory of Sobolev spaces}.
Contemp. Math., \textbf{338} (2003), 307–-340.



\bibitem{Kru} {\it Кругликов В.И.} Ёмкости конденсаторов и пространственные отображения,
квази\-кон\-формные в среднем // Матем. сб. -- 1986. --
\textbf{130}, № 2. -- C. 185-206.











\bibitem{Ge} {\it Gehring F.W.} Rings and quasiconformal
mappings in space // Trans. Amer. Math. Soc. -- 1962. --
\textbf{103}. -- P. 353--393.















































\bibitem{Fe} H. Federer, \emph{Geometric Measure Theory}, Springer,
Berlin etc., 1969.

\bibitem{Fu} {\it Fuglede B.} Extremal length and functional completion // Acta
Math. -- 1957. -- \textbf{98}. -- P. 171--219.





\bibitem{RR}
{\it Rado T., Reichelderfer P.V.} Continuous Transformations in Analysis. -- Springer--Verlag, Berlin, 1955. --
441 p.










\bibitem{Fe} {\it Федерер Г.} Геометрическая теория меры. -- Наука, Москва,
1987. -- 760 с.


















\bibitem{GoSa}{\it Golberg A.,  Salimov R.} Topological mappings of integrally
bounded p-moduli // Ann. Univ. Bucharest, Ser. Math. -- 2012. --
\textbf{3 (LXI)}, № 1. - P. 49-66.









\bibitem{Hes} {\it Hesse J.} A $p$-extremal length and $p$-capacity equality // Ark.
Mat. -- 1975. -- \textbf{13}. -- P. 131--144.








\bibitem{KR$_1$}
{\it Ковтонюк Д., Рязанов В.} К теории нижних $Q$-гомеоморфизмов //
Укр. мат. вiсник. -- 2008. -- \textbf{5}, № 2. -- С.~157--181.
























\bibitem{Re} {\it Решетняк Ю.Г.} Пространственные отображения с ограниченным
искажением. -- Наука, Новосибирск, 1982.






























\bibitem{Zi} {\it Ziemer W.P.} Extremal length and $p$-capacity. //  The Michigan Mathematical Journal 16 (1969), no. 1, 43--51.



\end{thebibliography}
\end{document}